\documentclass[10pt]{article}
\usepackage{graphicx}
\usepackage{amsmath}
\usepackage[dvips]{epsfig}
\usepackage{amssymb}

\makeatletter
\renewcommand{\section}{\@startsection
  {section}%
  {2}%
  {0mm}%
  {\baselineskip}%
  {0.5 \baselineskip}%
  {\centering}}
\makeatother
\begin{document}

\title { Note on $q$-extensions of  Euler numbers and polynomials of higher order}
\author{ Taekyun  Kim$^{1}$, Lee-Chae Jang$^{2}$, and Cheon-Seoung Ryoo$^{3}$  \\\\
$^1$EECS,
 Kyungpook National University, Taegu, 702-701,  Korea\\
         {\it   e-mail: tkim$@$knu.ac.kr } \\ \\
$^{2}$       Department of Mathematics  and Computer Science, \\
 KonKuk University, Chungju, Korea  \\
{\it  e-mail:   leechae-jang$@$hanmail.net}\\\\
$^{3}$         Department of Mathematics, \\
 Hannam University, Daejeon 306-791, Korea  \\
{\it  e-mail:   ryoocs$@$hnu.kr}
\\ \\
   }

\date{}
\maketitle

 {\footnotesize {\bf Abstract}\hspace{1mm}
In [14] Ozden-Simsek-Cangul constructed generating functions of
higher-order twisted $(h,q)$-extension of Euler polynomials and
numbers, by using $p$-adic $q$-deformed fermionic integral on $\Bbb
Z_p$. By applying their generating functions, they derived the
complete sums of products of the twisted $(h,q)$-extension of Euler
polynomials and numbers, see[13, 14]. In this paper we cosider the
new $q$-extension of Euler numbers and polynomials to be different
which is treated by Ozden-Simsek-Cangul. From our $q$-Euler numbers
and polynomials we derive some interesting identities and we
construct $q$-Euler zeta functions which interpolate the  new
$q$-Euler numbers and polynomials at a negative integer. Furthermore
we study Barnes' type $q$-Euler zeta functions. Finally we will
derive the new formula for " sums products of $q$-Euler numbers and
polynomials" by using fermionic $p$-adic $q$-integral on $\Bbb Z_p$.
 }

\medskip


\medskip


\bigskip

\begin{center}
 \bf{ \Large{1. Introduction and notations}}
 \end{center}

\bigskip
Throughout this paper we use the following notations. By
$\mathbb{Z}_p$ we  denote the ring of $p$-adic rational integers,
 $\Bbb Q$ denotes the field of rational numbers,
   $\mathbb{Q}_p$ denotes the field of $p$-adic rational numbers,  $\mathbb{C}$ denotes
  the complex number field, and $\mathbb{C}_p$ denotes the
  completion of algebraic closure of   $\mathbb{Q}_p$.
Let $\nu_p$ be the normalized exponential valuation of
$\mathbb{C}_p$ with $|p|_p=p^{-\nu_p(p)}=p^{-1}.$ When one talks of
$q$-extension, $q$ is considered in many ways such as an
indeterminate, a complex number $q\in  \mathbb{C},$ or $p$-adic
number $q\in\Bbb C_p .$ If $q\in \Bbb C$ one normally assume that
$|q|<1.$  If $q\in  \mathbb{C}_p,$ we normally assume that
$|q-1|_p<p^{-\frac{1}{p-1}}$ so that $q^x=\exp(x\log q)$ for
$|x|_p\leq 1.$ In this paper we use the following notation:
$$[x]_q =[x:q]=\frac{1-q^x}{1-q}, \mbox{ cf. [1, 2, 3, 4, 5, 21] }.$$
Hence,  $\lim_{q\rightarrow1}[x]=x$ for any $x$ with $|x|_p \leq 1
$ in the present $p$-adic case. Let $d$ be a fixed integer and let
$p$ be a fixed prime number. For any positive integer $N$, we set
$$\aligned
&X=\varprojlim_N ( \mathbb{Z}/dp^N \mathbb{Z}),\\
&X^*=\bigcup_{\begin{subarray}{l} 0<a<dp \\ (a,p)=1 \end{subarray}} ( a+dp \mathbb{Z}_p),\\
&a+dp^N \mathbb{Z}_p=\{x\in X\mid x\equiv a\pmod{dp^N}\},
\endaligned$$
where $a\in \mathbb{ Z}$ lies in $0\leq a<dp^N$. For any positive
integer $N$,
$$ \mu_q(a+dp^N \mathbb{Z}_p)=\dfrac{q^a}{[dp^N]_q}$$
is known to be a distribution on $X$, cf.[1-19]. From this
distribution we derive the $p$-adic $q$-integral on $\Bbb Z_p$ as
follows:
$$\int_{\Bbb Z_p} f(x)d\mu_q(x)=\lim_{N\rightarrow
\infty}\frac{1}{[p^N]_q} \sum_{x=0}^{p^N-1}q^xf(x), \text{ where
 $f\in UD(\Bbb Z_p)$ }, \text{ see [1-22] }.$$
Higher-order twisted Bernoulli and Euler numbers and polynomials are
studied by many authors ( see for detail [1-20]). In [14]
Ozden-Simsek-Cangul constructed generating functions of higher-order
twisted $(h,q)$-extension of Euler polynomials and numbers, by using
$p$-adic $q$-deformed fermionic integral on $\Bbb Z_p$. By applying
their generating functions, they derived the complete sums of
products of the twisted $(h,q)$-extension of Euler polynomials and
numbers, see[13, 14]. In this paper we cosider the new $q$-extension
of Euler numbers and polynomials to be different which is treated by
Ozden-Simsek-Cangul. From our $q$-Euler numbers and polynomials we
derive some interesting identities and we construct $q$-Euler zeta
functions which interpolate the  new $q$-Euler numbers and
polynomials at a negative integer. Furthermore we study Barnes' type
$q$-Euler zeta functions. Finally we will derive the new formula for
" sums products of $q$-Euler numbers and polynomials" by using
fermionic $p$-adic $q$-integral on $\Bbb Z_p$.

\medskip

\begin{center}
 \bf{ \Large{2. $q$-Extension of Euler numbers}}
 \end{center}
\medskip

In this section we assume that $q \in \mathbb{C}$ with $|q|<1$.
Now we consider the  $q$-Extension of Euler polynomials as
follows:
$$  F_{q}(x,t)=\dfrac{[2]_q }{qe^t+1} e^{xt} =\sum_{n=0}^{\infty}  \dfrac{ E_{n, q
}(x)}{n!} {t^n}, \quad  |t +\log q |< \pi.  \eqno(2.1)
$$
Note that $$ \lim_{q \to 1} F_q(x,t)=F(x,t)= \dfrac{2 }{e^t+1}
e^{xt} =\sum_{n=0}^{\infty}  \dfrac{ E_{n }(x)}{n!} {t^n}.$$ In
the special case $x=0$ the $q$-Euler polynomial
$E_{n,q}(0)=E_{n,q}$ will be called $q$-Euler numbers. It is easy
to see that $F_{q}(x,t)$ is analytic function in $\mathbb{C}$.
Hence we have
$$\sum_{n=0}^{\infty}  \dfrac{ E_{n, q
}(x)}{n!} {t^n}=\dfrac{[2]_q }{qe^t+1} e^{xt} =[2]_q
\sum_{n=0}^\infty (-1)^n q^n e^{(n+x)t}. \eqno(2.2)$$ If we take
the $k-$th derivative at $t=0$ on bosides in (2.2), then we have
$$E_{k,q}(x)=[2]_q \sum_{n=0}^\infty (-1)^n q^n (n+x)^k. \eqno(2.3)$$
From (2.3) we can define $q$-zeta function  which interpolating
$q$-Euler numbers at negative integer as follows:

For $s \in \mathbb{C}, $ we define
$$\zeta_{q}(s,x)=[2]_q \sum_{n=0}^\infty \dfrac{(-1)^n q^n}{(n+x)^s}, \quad s \in \mathbb{C}. \eqno(2.4)$$
Note that $\zeta_q(s,x)$ is analytic in complex $s$-plane. If we
take $s=-k( k \in \mathbb{Z}_+),$ then we have
$\zeta_q(-k,x)=E_{k,q}(x).$

By (2.3) and (2.4), we obtain  the following:

\bigskip

{\bf Theorem 1}.  For $k\in \mathbb{Z}_+$, we have
$$
E_{k, q}(x) =   [2]_q \sum_{n=0}^\infty {(-1)^n q^n}{(n+x)^k}.
$$

\medskip
Let $F_q(0,t)=F_q(t).$ Then we see that

$$ \aligned
  {[2]_q} \sum_{k=0}^{n-1}(-1)^k q^k e^{kt} &=
\dfrac{[2]_q}{1+qe^t}- [2]_q \dfrac{(-1)^n q^n
e^{nt}}{1+qe^t}\\
&=F_q(t)-(-1)^n q^n F_q(n,t).
 \endaligned \eqno(2.5)
$$
From (2.5), we derive
$$\sum_{k=0}^\infty \left( [2]_q  \sum_{k=0}^{n-1} (-1)^l q^l l^k \right)  \dfrac{ t^k }{k!}
=\sum_{k=0}^\infty \left( E_{k,q}-(-1)^n q^n E_{k,q}(n)  \right)
\dfrac{ t^k }{k!} . \eqno(2.6)$$ By comparing the coefficients on
both sides in (2.6),  we obtain the following:

\bigskip

{\bf Theorem 2}.  Let $n\in \Bbb N, k \in \mathbb{Z}_+$.  If $n
\equiv 0$ (mod 2), then we have $$ E_{k,q}-q^n E_{k,q}(n)=[2]_q
\sum_{k=0}^{n-1} (-1)^l q^l l^k .$$  If $n \equiv 1 $ (mod 2),
then we have $$ E_{k,q}+q^n E_{k,q}(n)=[2]_q \sum_{k=0}^{n-1}
(-1)^l q^l l^k.$$

\medskip

For $w_1, w_2, \cdots, w_r \in \mathbb{C},$ we consider the
multiple $q$-Euler polynomials of Barnes type as follows:
$$ \aligned
& F_q^r(w_1, w_2, \cdots, w_r |x, t)  = \dfrac{[2]_q^r
e^{xt}}{(qe^{w_1t}+1)(qe^{w_2t}+1)\cdots (qe^{w_rt}+1)} \\
& = \sum_{n=0}^\infty E_{n,q}(w_1, \cdots, w_r|x) \dfrac{t^n}{n!},
 \text{ where } \max_{ 1 \leq i \leq r} | w_i t+ \log q | < \pi
\endaligned \eqno(2.7)$$
For $x=0, E_{n,q}(w_1, \cdots, w_r|0)=E_{n,q}(w_1, \cdots, w_r)$
will be called the multiple $q$-Euler numbers of Barnes type. It
is easy to see that $F_q^r(w_1, w_2, \cdots, w_r |x, t)$ is
analytic function in the given region.  From (2.7) we derive
$$ [2]_q^r \sum_{n_1, \cdots, n_r=0}^\infty (-q)^{\sum_{i=1}^r n_i}e^{( \sum_{i=1}^r n_i w_i+x)t}=
\sum_{n=0}^\infty E_{n,q}(w_1, \cdots , w_r|x) \dfrac{t^n}{n!}.
\eqno(2.8)$$ By the $k$-th differentiation on both sides in (2.8),
we see that
$$[2]_q^r \sum_{n_1, \cdots, n_r=0}^\infty (-q)^{\sum_{i=1}^r n_i} \left(  \sum_{i=1}^r n_i w_i+x
\right)^k = E_{k,q}(w_1, \cdots , w_r|x) . \eqno(2.9)$$ From
(2.8), we can derive the following Barnes' type multiple $q$-Euler
zeta function as follows:

For $s \in \mathbb{C}$, define
$$ \zeta_{r,q}(w_1, w_2, \cdots, w_r | s,x)
=\sum_{n_1, \cdots, n_r=0}^\infty \dfrac{(-1)^{n_1+\cdots+n_r}
q^{n_1+\cdots+n_r} }{(n_1w_1+n_2w_2+\cdots + n_rw_r+x)^s}.
\eqno(2.10)$$ By (2.9) and (2.10), we obtain the following:

\bigskip

{\bf Theorem 3}.  For $k \in \mathbb{Z}_+,  w_1, w_2, \cdots, w_r
\in \mathbb{C}, $ we have  $$ \zeta_{r,q}(w_1, w_2, \cdots, w_r |
-k,x)=E_{k,q}(w_1, w_2, \cdots, w_r | x). $$

\medskip

Let $\chi$ be the primitive Drichlet character with conduct
$f$(=odd) $\in \mathbb{N}$. Then we consider a generalized Euler
numbers attached to $\chi$ as follows:
$$F_{\chi, q}(t)=\dfrac{[2]_q \sum_{a=0}^{f-1}  (-1)^a q^a \chi(a)
 e^{at}}{q^f e^{ft} +1}=\sum_{n=0}^\infty E_{n,\chi, q}
\dfrac{t^n}{n!},\eqno(2.11)$$  where $ | \log q+t|<
\frac{\pi}{f}$. The numbers $E_{n, \chi, q}$ will be called the
generalized $q$-Euler numbers attached to $\chi$. From (2.11), we
note that
$$ \aligned
& F_{\chi, q}( t)  = \dfrac{[2]_q \sum_{a=0}^{f-1}  (-1)^a q^a
\chi(a) e^{at}}{q^f e^{ft} +1}\\
& = [2]_q\sum_{a=0}^{f-1}  (-1)^a q^a \chi(a) \sum_{n=0}^\infty
q^{nf} (-1)^n e^{(a+nf)t}\\
&=[2]_q \sum_{n=0}^\infty \sum_{a=0}^{f-1}  (-1)^{a+nf} q^{a+nf} \chi(a+nf) e^{(a+nf)t}\\
&=[2]_q \sum_{n=0}^\infty (-1)^n q^n \chi(n) e^{nt}=
\sum_{n=0}^\infty E_{n,\chi, q} \dfrac{t^n}{n!} .
\endaligned
$$
Thus, we have
$$
 E_{n,\chi, q}=\left. \frac{d^k}{ dt^k} F_{\chi,
q}(t)\right\vert_{t=0} = [2]_q \sum_{n=1}^\infty (-1)^n q^n
\chi(n) n^k,  (k \in \mathbb{N}).
$$
Therefore , we can define  the   Dirichlet's type $l$-function
which interpolates at negative integer as follows:

For $ s \in \mathbb{C}$, we define $l_q(s,\chi)$ as
$$l_q(s,\chi)=[2]_q \sum_{n=1}^\infty \dfrac{(-1)^n q^n \chi(n)}{n^s}.
\eqno(2.13)$$ By (2.12) and (2.13), we obtain the following:

\bigskip

{\bf Theorem 4}.  For $k \in \mathbb{Z}_+, $ we have
$$ l_q(-k,\chi) = E_{k,\chi, q}. $$

\medskip

From (2.1) and the definition of $q$-Euler numbers, we derive
$$ \aligned
& F_{q}( t, x)  = \dfrac{[2]_q }{q e^{t} +1}e^{xt}  =
\sum_{n=0}^{\infty}E_{n,q} \dfrac{t^n}{n!}  \sum_{l=0}^\infty
 \dfrac{x^l}{l!} t^l \\
&= \sum_{m=0}^\infty \left( \sum_{n=0}^m E_{n,q} \binom mn x^{n-l}
\right) \dfrac{t^m}{m!}.
\endaligned  \eqno(2.14)
$$
By (2.14) we see that
$$E_{n,q}(x)= \sum_{m=0}^n E_{m,q} \binom nm x^{n-m}, \quad n \in \mathbb{Z}_+. \eqno(2.15)$$
For $f$ (=odd) $\in \mathbb{N}$, we note that
$$ \aligned
\sum_{n=0}^\infty E_{n,\chi, q} \dfrac{t^n}{n!} &  =
\dfrac{[2]_q }{q e^{t} +1}e^{xt} \\
&= [2]_q  \dfrac{1}{q^f e^{ft}+1 }\sum_{a=0}^{f-1}  (-1)^a q^a e^{( \dfrac{a+x}{f}) ft}  \\
& = \dfrac{[2]_q}{[2]_{q^f}}  \sum_{a=0}^{f-1}  (-1)^a q^a  \left(  \dfrac{ [2]_{q^f} e^{ ( \dfrac{a+x}{f} ) ft }} { q^f e^{ft} +1} \right) \\
&= \dfrac{[2]_q}{[2]_{q^f}} \sum_{a=0}^{f-1}  (-1)^a q^a
\sum_{n=0}^\infty E_{n, q^f} \left( \dfrac{a+x}{f}  \right)
\dfrac{ f^n t^n}{n!} .
\endaligned
$$
Thus, we have the distribution relation for $q$-Euler polynomials
as follows:

\bigskip

{\bf Theorem 5}.  For $f$(=odd) $ \in \mathbb{N}, $ we have
$$  E_{n, q}(x)=\dfrac{f^n [2]_q}{[2]_{q^f}} \sum_{a=0}^{f-1}  (-1)^a q^a
 E_{n, q^f} \left( \dfrac{a+x}{f}  \right). $$

\medskip

For $k, n \in \mathbb{N}$ with $n \equiv 0$(mod 2), it is easy to
see that
$$ \aligned
& [2]_q \sum_{l=0}^{n-1}(-1)^{l-1} q^l l ^k  = q^n
E_{k,q}(n)-E_{k,q} \\
&= q^n  \sum_{m=0}^k \binom km n^{k-m} E_{m,q}-E_{k,q} = q^n
\sum_{m=0}^{k-1} \binom km  E_{m,q} n^{k-m} + (q^n-1)E_{k,q}.
\endaligned
$$
Therefore, we obtain the following

\bigskip

{\bf Theorem 6}. For $k, n \in \mathbb{N}$ with $n \equiv 0$(mod
2),  we have
$$  [2]_q \sum_{l=0}^{n-1}(-1)^{l-1} q^l l ^k  = q^n
\sum_{m=0}^{k-1} \binom km  E_{m,q} n^{k-m} + (q^n-1)E_{k,q}. $$

\medskip

\begin{center}
 \bf{ \Large{3. Witt's type formulae on $\mathbb{Z}_p$ in $p$-adic number field}}
 \end{center}

\medskip

In this section we assume that  $q\in\Bbb C_p $ with $|1-q|_p<1$.
$g$ is a  uniformly differentiable function  at  a point $a \in
\mathbb{Z}_p $, and write  $g \in UD(\mathbb{Z}_p)$, if the
difference quotient
$$F_g(x,y)=\dfrac{g(x)-g(y)}{x-y},$$ has a limit $f'(a)$ as $(x,y)
\to (a,a)$.  For $g \in UD(\mathbb{Z}_p)$, an invariant $p$-adic
$q$-integral is defined as
$$I_{q}(g)=\int_{\mathbb{Z}_p}g(x) d\mu_{q}(x)=
\lim_{N \to \infty} \dfrac{1}{[p^N]_q }\sum_{x= 0}^{ p^N-1
}g(x)q^x.$$ The fermionic $p$-adic $q$-integral is also defined as
$$I_{-q}(g)=\int_{\mathbb{Z}_p}g(x) d\mu_{-q}(x)=
\lim_{N \to \infty}\dfrac{[2]_q }{1+q^{p^N} } \sum_{x= 0}^{ p^N-1
}g(x)(-1)^x q^x, \text{ see [4] }. \eqno(3.1)$$ From (3.1), we have
the integral equation as follows:
$$ q I_{-q}(g_1)+I_{-q}(g)=[2]_q f(0),  \text{ where } g_1(x)=g(x+1). \eqno(3.2)$$
If we take $f(x)=e^{tx}$, then we have
$$I_q(e^{tx})= \int_{\mathbb{Z}_p}e^{x t } d\mu_{-q}(x)=\dfrac{[2]_q}{qe^t+1}. \eqno(3.3)$$
From (3.3), we note that
$$ \sum_{n=0}^\infty \int_{\mathbb{Z}_p}x^n d\mu_{-q}(x) \dfrac{t^n}{n!} =\dfrac{[2]_q}{qe^t+1}= \sum_{n=0}^\infty E_{n,q} \dfrac{t^n}{n!}. $$
By comparing the coefficient on the both sides, we see that
$$  \int_{\mathbb{Z}_p}x^n d\mu_{-q}(x) = E_{n,q}, n \in \mathbb{Z}_+. \eqno(3.4) $$
By the same method we see that
$$  \int_{\mathbb{Z}_p} e^{(x+y)t}  d\mu_{-q}(x) = \dfrac{[2]_q}{qe^t+1} e^{xt} =\sum_{n=0}^\infty E_{n,q}(x) \dfrac{t^n}{n!}. $$
Hence, we have the formula of Witt's type for $q$-Euler polynomial
as follows:
$$  \int_{\mathbb{Z}_p}(x+y)^n d\mu_{-q}(x) = E_{n,q}, n \in \mathbb{Z}_+.  $$
For $n \in \mathbb{Z}_+$, let $g_n(x)=g(x+n)$.  Then we have
$$ q^n I_{-q}(g_n)+(-1)^{n-1}I_{-q}(g)=[2]_q   \sum_{l=0}^{n-1} (-1)^{n-1-l} q^l g(l). \eqno(3.5)$$
If $n$ is odd positive integer, then we have
$$ q^n I_{-q}(g_n)+ I_{-q}(g)=[2]_q   \sum_{l=0}^{n-1} (-1)^{l} q^l g(l). \eqno(3.6)$$

Let $\chi$ be the primitive Drichlet character with conduct
$f$(=odd) $\in \mathbb{N}$ and let $g(x)=\chi(x) e^{xt} $. From
(3.3) we derive
$$I_{-q}({\chi(x) e^{xt}})=\dfrac{[2]_q \sum_{a=0}^{f-1}  (-1)^a q^a \chi(a)
 e^{at}}{q^f e^{ft} +1}=\sum_{n=0}^\infty E_{n,\chi, q}
\dfrac{t^n}{n!}.$$ Thus, we have the Witt's formula for
generalized $q$-Euler numbers attached to $\chi$ as follows:
$$  \int_{\mathbb{Z}_p} \chi(x) x^n d\mu_{-q}(x) = E_{n,\chi, q}, n  \geq 0.  $$

\medskip

\begin{center}
 \bf{ \Large{4. Higher order  $q$-Euler numbers and polynomials}}
 \end{center}
\medskip

In this section we also assume that $q\in\Bbb C_p $ with
$|1-q|_p<1$. Now we study on higher order $q$-Euler numbers and
polynomials and sums of products of $q$-Euler numbers. First, we try
to consider the multivariate  fermionic $p$-adic $q$-integral on $
\mathbb{Z}_p $ as follows:

$$ \aligned
& \underset{r \mbox{ times}}{\underbrace{ \int_{ \mathbb{Z}_p}
\cdots \int_{\mathbb{Z}_p }}} e^{ (a_1x_1+a_2x_2+ \cdots + a_r
x_r)t} e^{xt} d\mu_{-q}(x_1) \cdots d\mu_{-q}(x_r) \\
&= \dfrac{[2]_q^r}{(qe^{a_1t}+1)(q e^{a_2t}+1) \cdots (qe^{a_r
t}+1)} e^{xt},
\endaligned
\eqno(4.1) $$  where $ a_1, a_2, \cdots, a_r \in \mathbb{Z}_p$.

From (4.1) we consider the multiple $q$-Euler polynomials as
follows:
$$ \dfrac{[2]_q^r}{(qe^{a_1t}+1)(q e^{a_2t}+1) \cdots (qe^{a_r
t}+1)} e^{xt}= \sum_{n=0}^\infty E_{n,q}( a_1, a_2, \cdots, a_r
|x) \dfrac{t^n}{n!}. \eqno(4.2) $$ In the special case $( a_1,
a_2, \cdots, a_r)=(1, 1, \cdots, 1),$ we write
$$E_{n,q} ( \underset{r \mbox{ times}}{\underbrace{ a_1
\cdots a_r}}|x )=E_{n,q}^{(r)}(x).$$ For $x=0$, the multiple
$q$-Euler polynomials will be called as $q$-Euler numbers of order
$r$.

For (4.2) and (4.2) we note that
$$ E_{n,q}( a_1, a_2, \cdots, a_r
|x)=\underset{r \mbox{ times}}{\underbrace{ \int_{ \mathbb{Z}_p}
\cdots \int_{\mathbb{Z}_p }}}(a_1 x_1+ \cdots + a_r x_r+x)^n
\prod_{j=1}^r d \mu_{-q} (x_r).$$ It is easy to check that
$$E_{n,q}( a_1, a_2, \cdots, a_r
|x)=\sum_{l=0}^n \binom nl x^{n-l} E_{l,q}( a_1, a_2, \cdots,
a_r),$$ where $ E_{n,q}( a_1, a_2, \cdots, a_r )= E_{n,q}( a_1,
a_2, \cdots, a_r |0)$.  Multinomial  theorem is well known as
follows:
$$
\left( \sum_{j=1}^r x_j \right)^n = \sum_{\substack
{l_1,\cdots,l_r \geq 0 \\ l_1+ \cdots +l_r=n}}  \binom{n}{l_1,
\cdots, l_r} \prod_{a=1}^r x_a^{l_a}, $$  where $$\binom
{n}{l_1,\cdots,l_r}=\dfrac{n!}{l_1! l_2! \cdots l_r !}.$$ By (4.2)
and (4.3) we easily see that
$$E_{n,q}^{(r)}(x)= \sum_{m=0}^n \sum_{\substack{l_1,\cdots,l_r
\geq 0 \\  l_1+ \cdots +l_r=m}}  \binom nm \binom
{m}{l_1,\cdots,l_r} x^{n-m}\prod_{j=1}^r E_{l_j,q}. $$


\bigskip


\begin{thebibliography}{99}

\bibitem{}
 { M. Cenkci},
{\it The $p$-adic generalized twisted $(h,q)$-Euler-$l$-function and
its applications},
  { Advan. Stud. Contemp. Math.}, 15(2007), 37-47.

\bibitem{}
 { L.C. Jang, S.D. Kim, D.-W. Park, Y.-S. Ro},
{\it  A note on Euler Numbers and Polynomials},
  {J. Ineq. Appl.}, Vol 2006(Art.ID 34602)(2006), 1-5.

\bibitem{}
 { T. Kim},
{\it A note on $p$-adic invariant integral in the rings of $p$-adic
integer },
  { Advan. Stud. Contemp. Math.}, 13(2006), 95-99.


\bibitem{}
 { T. Kim},
{\it A Note on $p$-Adic $q$-integral on $\Bbb Z_p$ Associated with
$q$-Euler Numbers },
  { Advan. Stud. Contemp. Math.}, 15(2007), 133-137.


\bibitem{}
 { T. Kim, M.S. Kim, L.C. Jang, S. H. Rim},
{\it New $q$-Euler Numbers and Polynomials Associated with $p$-Adic
$q$-Integrals},
  { Advan. Stud. Contemp. Math.}, 15(2007), 243-252.

\bibitem{}
 { T. Kim},
{\it A new approach to $q$-zeta function},
  { Advan. Stud. Contemp. Math.}, 11(2005), 157-162.

  \bibitem{}
 { T. Kim},
{\it A note on exploring the sums of powers of consequtive
$q$-integer},
  { Advan. Stud. Contemp. Math.}, 11(2005), 137-140.

  \bibitem{}
 { T. Kim},
{\it Sums of powers of consequtive $q$-integer},
  { Advan. Stud. Contemp. Math.}, 9(2004), 15-18.


\bibitem{}
 { T. Kim},
{\it $q$-generalized Euler numbers and polynomials },
  { Russ. J. Math. Phys.}, 13 (2006), 293-298.

\bibitem{}
 { T. Kim},
{\it Multiple $p$-adic $L$-function},
  { Russ. J. Math. Phys.}, 13 (2006), 151-157.


\bibitem{}
 { T. Kim},
{\it $q$-Volkenborn integration},
  { Russ. J. Math. Phys.},  9 (2002), 288-298.


\bibitem{}
 { T. Kim},
{\it  A note on $q$-Volkenborn integration},
  {Proc. Jangjeon Math. Soc.}, Vol 2(2001), 45-49.


\bibitem{}
{ H. Ozden, Y. Simsek, S.H. Rim, I. Cangul}, {\it A note on $p$-adic
$q$-Euler measure}, {Advan. Stud. Contemp. Math.}, 14(2007),
233-239.

\bibitem{}
 { H. Ozden, Y. Simsek, I.N. Cangul},
{\it  Remarks on sum of products of $(h,q)$-twisted Euler
polynomials and numbers},
  {J. Ineq. Appl.(2007), in press: http://www.hindawi.com/journals/jia/oct.2007.html}.


\bibitem{}
{ Y. Simsek}, {\it  On $p$-adic twisted $q$-$L$-function related to
generalized twisted Bernoulli numbers}, { Russ. J. Math. Phys.}, 13
(2006), 340-348.


\bibitem{}
{ Y. Simsek}, {\it  Theorem  on twisted $L$-function  and twisted
 Bernoulli numbers}, { Advan. Stud. Contemp. Math.}, 11 (2005), 205-258.

\bibitem{}
{ Y. Simsek}, {\it  Twisted $(h,q)$-Bernoulli numbers and
polynomials related to twisted $(h,q)$-zeta function and
$L$-function}, { J. Math. Anal. Appl.}, 324 (2006), 790-804.


\bibitem{}
{ C. S. Ryoo,   H. Song, ,   R. P. Agarwal},  {\it On the roots of
the $q$-analogue of Euler-Barnes' polynomials},{ Advan. Stud.
Contemp. Math.}, 9(2004), 153-163.



\bibitem{}
{ C. S. Ryoo},  {\it A note on $q$-Bernoulli numbers and
polynomials }, { Appl. Math. Letters}, 20(2007), 524-531.

\bibitem{}
{ C. S. Ryoo},  {\it The zero of generalized twisted Bernoulli
polynomials }, { Adv. Theor. Appl. Math.}, 1(2006), 143-148.

\bibitem{}
{ C. S. Ryoo,   T. Kim,  R. P. Agarwal},  {\it Distribution of the
roots of the Euler-Barnes's type $q$-Euler polynomials}, { Neural
Parallel Sci. Comput.}, 13(2005), 381-392.

\bibitem{}
{ C. S. Ryoo,   T. Kim,  L.C. Jang},  {\it Some relationships
between the analogs of Euler numbers and polynomials}, {J. Ineq.
Appl.(2007), in press:
http://www.hindawi.com/journals/jia/aug.2007.html}.



\end{thebibliography}
\end{document}